\documentclass{amsart}
\usepackage{amsmath}
\usepackage{amsthm}
\usepackage{amsopn}
\usepackage{amsfonts}
\usepackage{amssymb}
\title{Splitting of Sharply 2-Transitive Groups of Characteristic 3}
\author{Seyfi T\"urkelli\\Istanbul Bilgi University\\Kustepe Istanbul Turkey}
\date{November 23, 2003.}
\newtheorem{thm}{Theorem}
\newtheorem{lemma}{Lemma}
\DeclareMathOperator{\chr}{char}

\begin{document}
\maketitle

\begin{abstract}
We give a group theoretic proof of the splitting of sharply
2-transitive groups of characteristic 3.
\end{abstract}\ \\


A \emph{sharply $2$-transitive group} is a pair $(G,X)$ where $G$
is a group acting on the set $X$ in such a way that for all
$x,\,y,\,z,\,t\in X$ such that $x\neq y$ and $z\neq t$ there is a
unique $g\in G$ for which $gx=z$ and $gy=t.$ From now on $(G,X)$
will stand for a sharply $2$-transitive group with $|X|\geq 3$. We
fix an element $x\in X$. We let $H:=\{g\in G: gx=x\}$ denote the
stabilizer of $x$. Finally we let $I$ denote the set of
involutions (elements of order 2) of $G$.

It follows easily from the definition that the group $G$ has an
involution, in fact any element of $G$ that sends a distinct pair
$(y,z)$ of $X$ to the pair $(z,y)$ is an involution by sharp
transitivity. It is also known that $I$ is one conjugacy class and
the nontrivial elements of $I^2$ cannot fix any point (See Lemma 1
and Lemma 4). Then one can see that $I^2$ cannot have an
involution if $H$ has an involution.

In case $H$ has no involution, one says that $\chr(G)=2$.

Let us assume that $\chr(G)\neq 2$. Then $I^2\setminus\{1\}$ is
one conjugacy class \cite[Lemma 11.45]{bn}. Since $I^2$ is closed
under power taking, either the nontrivial elements of $I^2$ all
have order $p$ for some prime $p\neq 2$ or $I^2$ has no nontrivial
torsion element. One writes $\chr(G)=p$ or $\chr(G)=0$ depending
on the case.

One says that $G$ \emph{splits} if the one point stabilizer $H$
has a normal complement in $G$. It is not known whether or not an
infinite sharply 2-transitive group splits, except for those of
characteristic 3. Results in this direction for some special cases
can be found in \cite[\S 11.4]{bn} and \cite[ch 2]{ke}. We will
prove that if $\chr(G)=3$ then $G$ splits, a result of W.\ Kerby
\cite[Theorem 8.7]{ke}. But Kerby's proof is in the language of
near domains and is not easily accessible. Here, we give a much
simpler proof of this fact, in fact an experienced reader can
directly go to the proof the Theorem, which contains only a simple
computation (all the lemmas are well-known facts).

All the results of this short and elementary paper can be found in
\cite[\S 11.4]{bn}, except for the final theorem.

\begin{lemma}
$I$ is one conjugacy class.
\end{lemma}
\begin{proof}
Let $i,j\in I$ and $x\in X$ be such that $jx\neq x$ and $ix\neq
x$. Since $G$ is $2-$transitive, there exists a $g\in G$ such that
$gx=x$ and $gjx=ix$. Then $i^gjx=x$ and $i^gj(jx)=jx$. By double
sharpness of $G$, $i^gj=1$. Hence, $i^g=j$ and we are done.
\end{proof}

\begin{lemma}
    If $N$ is a nontrivial normal subgroup of
    $G$ then $G=NH$.
\end{lemma}
\begin{proof}
    Let $g\in G\setminus
    H, a\in N, y\in X\setminus \{x\}$ be such that $ay\neq y$ and
    $h\in G$ be such that $hx=y$ and $hgx=ay$. Then $(a^{-1})^hg\in
    H$ and $g\in NH$. Since $1\in N$, it holds for all $g\in G$.
\end{proof}
\begin{lemma}
    $H$ has at most one
    involution.
\end{lemma}
\begin{proof}
    Let $i,j\in H\cap I,\, y\in
    X\setminus \{x\},\, g\in G $ be such that $gjy=iy$ and $gy=y$. Then
    $ji^g(y)=y$ and $ji^g(jy)=jy$. Since $ji^g$ fixes two different
    points and $G$ is sharply $2$-transitive, $ji^g=1$ and $j=i^g$.
    One can easily see that $H\cap H^z\neq \{1\}$ if and only if $z\in H$.
    Therefore $g\in H$ as $j\in H\cap H^g$. Since $g$ fixes two
    points, namely $x$ and $y$, $g=1$. Hence $i=j$ and we are done.
\end{proof}
\begin{lemma}
    A nontrivial element of $I^2$ cannot fix any
    element of $X$.
\end{lemma}
\begin{proof}
    Assume not. Then, there are distinct involutions $i,\,j$ such that
    $ij$ fixes a point. Since $G$
    is transitive, we may assume $ij\in H$. It follows from Lemma 3 that $j\notin H$ otherwise $i\in H$,
    hence a contradiction.
    On the other hand, $(ij)^{-1}=(ji)=(ij)^j$ and $(ij)^j\in
    H\cap H^j$. Therefore, $j\in H$, a contradiction.
\end{proof}
\begin{lemma}
    If the elements of $Ii$ commute with each other for
    some $i\in I$, then $I^2$ is a normal subgroup of $G$.
\end{lemma}
\begin{proof}
    It suffices
    to prove that $I^2$ is closed under multiplication. Let
    $i,\,j,\,k,\,w\in I$. We claim that $ijkw\in I^2$. By Lemma 1,
    we may assume that the elements of $Ii$ commute with each other. Noting that $Ii = iI$, we have
    $(ijk)^2=ijkijk=kiijjk=1$. So, $ijk\in I \cup \{1\}$. If $ijk\in
    I$, we are done. Assume $ijk=1$. If $H$ has an involution, by
    Lemma 1,
    $(ij)^g=k^g\in H$ for some $g\in G$ , i.e.\ $(ij)^g$ fixes $x$, contradicting
    Lemma 4. If $H$ has no involution, $ij=k\in I$ and, by Lemma 1, $I\subseteq
    I^2$. Therefore, $ijkw=w\in I^2$.
\end{proof}
\begin{lemma}
    If $H$ has an involution, then the action of $G$ on
    $X$ is equivalent to the action of $G$ on $I$ by
    conjugation.
\end{lemma}
\begin{proof}
    Let $i\in H$ be an
    involution. It is easy to see that the action of $G$ on $X$ is
    equivalent to the action of $G$ on the left coset space $G/H$. So
    we may assume that the set $X$ is the left coset space $G/H$.
    Consider the map from $G/H$ to $I$ defined as $\bar{g}\mapsto i^{g^{-1}}$ for
    $g\in G$. One can easily see that
    this is the required equivalence.
\end{proof}
\begin{thm}
    If $\chr(G)=3$ then $G$ splits.
\end{thm}
\begin{proof}
    We claim that $G=I^2\rtimes H$. If $I^2$ is a normal subgroup of $G$,
    then we know that $H\cap I^2=\{1\}$ by Lemma 4 and $G=I^2H$ by Lemma
    2. Therefore, we just need to prove that $I^2$ is a normal subgroup of $G$.
    By lemma 5, it is enough to show that the elements of $Ii$ commute with each
    other for some $i\in I$. Let $i\in H\cap I$ be the (unique) involution of $H$ and let $ji,ki\in Ii$.
    We may assume that $j\neq k$.
    By double sharpness of $G$, it suffices to prove that $jiki$ and
    $kiji$ agree on two different points. By Lemma 6, we can take
    $X$ to be $I$ and the action to be the conjugation. We now claim
    that $jiki$ and $kiji$ agree on $j$ and $k$ i.e.\ that
    $j^{jiki}=j^{kiji}$ and $k^{jiki}=k^{kiji}$. By symmetry of the
    situation, it is enough to prove one of the equalities. Since
    $\chr(G)=3$, $i^j=j^i$ for all $i,j\in I$ and so we have
    $$j^{jiki}=j^{(k^i)}=(k^i)^j=k^{ij}=k^{jiji}=(k^j)^{iji}=(j^k)^{iji}=j^{kiji}.$$
\end{proof}

\bigskip
\ \\

\noindent Seyfi T\"URKELL\.{I}\\ Park Rheyngaerde 100 D 16\\ 3545
NE Utrecht The Netherlands\\ turkelli@wisc.edu

\end{document}